\catcode`\@=11

\font\tenmsa=msam10
\font\sevenmsa=msam7
\font\fivemsa=msam5
\font\tenmsb=msbm10
\font\sevenmsb=msbm7
\font\fivemsb=msbm5
\newfam\msafam
\newfam\msbfam
\textfont\msafam=\tenmsa  \scriptfont\msafam=\sevenmsa
  \scriptscriptfont\msafam=\fivemsa
\textfont\msbfam=\tenmsb  \scriptfont\msbfam=\sevenmsb
  \scriptscriptfont\msbfam=\fivemsb

\def\hexnumber@#1{\ifnum#1<10 \number#1\else
 \ifnum#1=10 A\else\ifnum#1=11 B\else\ifnum#1=12 C\else
 \ifnum#1=13 D\else\ifnum#1=14 E\else\ifnum#1=15 F\fi\fi\fi\fi\fi\fi\fi}

\def\msa@{\hexnumber@\msafam}
\def\msb@{\hexnumber@\msbfam}
\mathchardef\boxdot="2\msa@00
\mathchardef\boxplus="2\msa@01
\mathchardef\boxtimes="2\msa@02
\mathchardef\square="0\msa@03
\mathchardef\blacksquare="0\msa@04
\mathchardef\centerdot="2\msa@05
\mathchardef\lozenge="0\msa@06
\mathchardef\blacklozenge="0\msa@07
\mathchardef\circlearrowright="3\msa@08
\mathchardef\circlearrowleft="3\msa@09
\mathchardef\rightleftharpoons="3\msa@0A
\mathchardef\leftrightharpoons="3\msa@0B
\mathchardef\boxminus="2\msa@0C
\mathchardef\Vdash="3\msa@0D
\mathchardef\Vvdash="3\msa@0E
\mathchardef\vDash="3\msa@0F
\mathchardef\twoheadrightarrow="3\msa@10
\mathchardef\twoheadleftarrow="3\msa@11
\mathchardef\leftleftarrows="3\msa@12
\mathchardef\rightrightarrows="3\msa@13
\mathchardef\upuparrows="3\msa@14
\mathchardef\downdownarrows="3\msa@15
\mathchardef\upharpoonright="3\msa@16

\mathchardef\downharpoonright="3\msa@17
\mathchardef\upharpoonleft="3\msa@18
\mathchardef\downharpoonleft="3\msa@19
\mathchardef\rightarrowtail="3\msa@1A
\mathchardef\leftarrowtail="3\msa@1B
\mathchardef\leftrightarrows="3\msa@1C
\mathchardef\rightleftarrows="3\msa@1D
\mathchardef\Lsh="3\msa@1E
\mathchardef\Rsh="3\msa@1F
\mathchardef\rightsquigarrow="3\msa@20
\mathchardef\leftrightsquigarrow="3\msa@21
\mathchardef\looparrowleft="3\msa@22
\mathchardef\looparrowright="3\msa@23
\mathchardef\circeq="3\msa@24
\mathchardef\succsim="3\msa@25
\mathchardef\gtrsim="3\msa@26
\mathchardef\gtrapprox="3\msa@27
\mathchardef\multimap="3\msa@28
\mathchardef\therefore="3\msa@29
\mathchardef\because="3\msa@2A
\mathchardef\doteqdot="3\msa@2B

\mathchardef\triangleq="3\msa@2C
\mathchardef\precsim="3\msa@2D
\mathchardef\lesssim="3\msa@2E
\mathchardef\lessapprox="3\msa@2F
\mathchardef\eqslantless="3\msa@30
\mathchardef\eqslantgtr="3\msa@31
\mathchardef\curlyeqprec="3\msa@32
\mathchardef\curlyeqsucc="3\msa@33
\mathchardef\preccurlyeq="3\msa@34
\mathchardef\leqq="3\msa@35
\mathchardef\leqslant="3\msa@36
\mathchardef\lessgtr="3\msa@37
\mathchardef\backprime="0\msa@38
\mathchardef\risingdotseq="3\msa@3A
\mathchardef\fallingdotseq="3\msa@3B
\mathchardef\succcurlyeq="3\msa@3C
\mathchardef\geqq="3\msa@3D
\mathchardef\geqslant="3\msa@3E
\mathchardef\gtrless="3\msa@3F
\mathchardef\sqsubset="3\msa@40
\mathchardef\sqsupset="3\msa@41
\mathchardef\trianglerighteq="3\msa@44
\mathchardef\trianglelefteq="3\msa@45
\mathchardef\bigstar="0\msa@46
\mathchardef\between="3\msa@47
\mathchardef\blacktriangledown="0\msa@48
\mathchardef\blacktriangleright="3\msa@49
\mathchardef\blacktriangleleft="3\msa@4A
\mathchardef\blacktriangle="0\msa@4E
\mathchardef\triangledown="0\msa@4F
\mathchardef\eqcirc="3\msa@50
\mathchardef\lesseqgtr="3\msa@51
\mathchardef\gtreqless="3\msa@52
\mathchardef\lesseqqgtr="3\msa@53
\mathchardef\gtreqqless="3\msa@54
\mathchardef\Rrightarrow="3\msa@56
\mathchardef\Lleftarrow="3\msa@57
\mathchardef\veebar="2\msa@59
\mathchardef\barwedge="2\msa@5A
\mathchardef\doublebarwedge="2\msa@5B
\mathchardef\angle="0\msa@5C
\mathchardef\measuredangle="0\msa@5D
\mathchardef\sphericalangle="0\msa@5E
\mathchardef\varpropto="3\msa@5F
\mathchardef\smallsmile="3\msa@60
\mathchardef\smallfrown="3\msa@61
\mathchardef\Subset="3\msa@62
\mathchardef\Supset="3\msa@63
\mathchardef\Cup="2\msa@64

\mathchardef\Cap="2\msa@65

\mathchardef\curlywedge="2\msa@66
\mathchardef\curlyvee="2\msa@67
\mathchardef\leftthreetimes="2\msa@68
\mathchardef\rightthreetimes="2\msa@69
\mathchardef\subseteqq="3\msa@6A
\mathchardef\supseteqq="3\msa@6B
\mathchardef\bumpeq="3\msa@6C
\mathchardef\Bumpeq="3\msa@6D
\mathchardef\lll="3\msa@6E

\mathchardef\ggg="3\msa@6F

\mathchardef\circledS="0\msa@73
\mathchardef\pitchfork="3\msa@74
\mathchardef\dotplus="2\msa@75
\mathchardef\backsim="3\msa@76
\mathchardef\backsimeq="3\msa@77
\mathchardef\complement="0\msa@7B
\mathchardef\intercal="2\msa@7C
\mathchardef\circledcirc="2\msa@7D
\mathchardef\circledast="2\msa@7E
\mathchardef\circleddash="2\msa@7F
\def\ulcorner{\delimiter"4\msa@70\msa@70 }
\def\urcorner{\delimiter"5\msa@71\msa@71 }
\def\llcorner{\delimiter"4\msa@78\msa@78 }
\def\lrcorner{\delimiter"5\msa@79\msa@79 }
\def\yen{\mathhexbox\msa@55 }
\def\checkmark{\mathhexbox\msa@58 }
\def\circledR{\mathhexbox\msa@72 }
\def\maltese{\mathhexbox\msa@7A }
\mathchardef\lvertneqq="3\msb@00
\mathchardef\gvertneqq="3\msb@01
\mathchardef\nleq="3\msb@02
\mathchardef\ngeq="3\msb@03
\mathchardef\nless="3\msb@04
\mathchardef\ngtr="3\msb@05
\mathchardef\nprec="3\msb@06
\mathchardef\nsucc="3\msb@07
\mathchardef\lneqq="3\msb@08
\mathchardef\gneqq="3\msb@09
\mathchardef\nleqslant="3\msb@0A
\mathchardef\ngeqslant="3\msb@0B
\mathchardef\lneq="3\msb@0C
\mathchardef\gneq="3\msb@0D
\mathchardef\npreceq="3\msb@0E
\mathchardef\nsucceq="3\msb@0F
\mathchardef\precnsim="3\msb@10
\mathchardef\succnsim="3\msb@11
\mathchardef\lnsim="3\msb@12
\mathchardef\gnsim="3\msb@13
\mathchardef\nleqq="3\msb@14
\mathchardef\ngeqq="3\msb@15
\mathchardef\precneqq="3\msb@16
\mathchardef\succneqq="3\msb@17
\mathchardef\precnapprox="3\msb@18
\mathchardef\succnapprox="3\msb@19
\mathchardef\lnapprox="3\msb@1A
\mathchardef\gnapprox="3\msb@1B
\mathchardef\nsim="3\msb@1C
\mathchardef\napprox="3\msb@1D
\mathchardef\nsubseteqq="3\msb@22
\mathchardef\nsupseteqq="3\msb@23
\mathchardef\subsetneqq="3\msb@24
\mathchardef\supsetneqq="3\msb@25
\mathchardef\subsetneq="3\msb@28
\mathchardef\supsetneq="3\msb@29
\mathchardef\nsubseteq="3\msb@2A
\mathchardef\nsupseteq="3\msb@2B
\mathchardef\nparallel="3\msb@2C
\mathchardef\nmid="3\msb@2D
\mathchardef\nshortmid="3\msb@2E
\mathchardef\nshortparallel="3\msb@2F
\mathchardef\nvdash="3\msb@30
\mathchardef\nVdash="3\msb@31
\mathchardef\nvDash="3\msb@32
\mathchardef\nVDash="3\msb@33
\mathchardef\ntrianglerighteq="3\msb@34
\mathchardef\ntrianglelefteq="3\msb@35
\mathchardef\ntriangleleft="3\msb@36
\mathchardef\ntriangleright="3\msb@37
\mathchardef\nleftarrow="3\msb@38
\mathchardef\nrightarrow="3\msb@39
\mathchardef\nLeftarrow="3\msb@3A
\mathchardef\nRightarrow="3\msb@3B
\mathchardef\nLeftrightarrow="3\msb@3C
\mathchardef\nleftrightarrow="3\msb@3D
\mathchardef\divideontimes="2\msb@3E
\mathchardef\varnothing="0\msb@3F
\mathchardef\nexists="0\msb@40
\mathchardef\mho="0\msb@66
\mathchardef\thorn="0\msb@67
\mathchardef\beth="0\msb@69
\mathchardef\gimel="0\msb@6A
\mathchardef\daleth="0\msb@6B
\mathchardef\lessdot="3\msb@6C
\mathchardef\gtrdot="3\msb@6D
\mathchardef\ltimes="2\msb@6E
\mathchardef\rtimes="2\msb@6F
\mathchardef\shortmid="3\msb@70
\mathchardef\shortparallel="3\msb@71
\mathchardef\smallsetminus="2\msb@72
\mathchardef\thicksim="3\msb@73
\mathchardef\thickapprox="3\msb@74
\mathchardef\approxeq="3\msb@75
\mathchardef\succapprox="3\msb@76
\mathchardef\precapprox="3\msb@77
\mathchardef\curvearrowleft="3\msb@78
\mathchardef\curvearrowright="3\msb@79
\mathchardef\digamma="0\msb@7A
\mathchardef\varkappa="0\msb@7B
\mathchardef\hslash="0\msb@7D
\mathchardef\hbar="0\msb@7E
\mathchardef\backepsilon="3\msb@7F
\def\Bbb{\ifmmode\let\next\Bbb@\else
 \def\next{\errmessage{Use \string\Bbb\space only in math mode}}\fi\next}
\def\Bbb@#1{{\Bbb@@{#1}}}
\def\Bbb@@#1{\fam\msbfam#1}

\catcode`\@=\active

\def\CF{\hbox{{$\cal F$}}}

\def\CM{\hbox{{$\cal M$}}}
\def\CN{\hbox{{$\cal N$}}}

\def\R{{\Bbb R}}
\def\C{{\Bbb C}}
\def\Z{{\Bbb Z}}

\def\lform{\hbox{$\sqcup$}\llap{\hbox{$\sqcap$}}}

\def\h{{{1\over2}}}
\def\eps{{\epsilon}}

\def\<{\langle}
\def\>{\rangle}

\def\tens{\mathop{\otimes}}

\def\isom{{\cong}}

\def\Ad{{\rm Ad}}

\def\id{{\rm id}}

\def\o{{}_{(1)}}\def\t{{}_{(2)}}

\def\extd{{\rm d}}

\def\proof{\goodbreak\noindent{\bf Proof\quad}}
\def\endproof{{\ $\lform$}\bigskip }

\def\text#1{{\rm #1}}
\def\note#1{}

\def\nquad{{\!\!\!\!\!\!}}

\def\equad{\nquad}
\def\eqn#1#2{\begin{equation}#2\label{#1}\end{equation}}

\def\align#1{\begin{eqnarray*}#1\end{eqnarray*}}

\def\cmath#1{\[\begin{array}{c} #1 \end{array}\]}

\def\ceqn#1#2{\begin{equation}\label{#1}\begin{array}{c}#2\end{array}
\end{equation}}

\documentstyle[11pt]{article}
\newtheorem{lemma}{Lemma}[section]
\newtheorem{propos}[lemma]{Proposition}

\begin{document}
\baselineskip 23pt

{\ }\hskip 4.7in   Damtp/97-62
\vspace{.2in}

\begin{center} {\Large QUANTUM GEOMETRY OF
FIELD EXTENSIONS}
\\ \baselineskip 13pt{\ }
{\ }\\ Shahn Majid\footnote{Royal Society University Research
Fellow and Fellow of Pembroke College, Cambridge, England.}\\ {\
}\\ Department of Applied Mathematics and Theoretical Physics\\
University of Cambridge, Cambridge CB3 9EW, UK\\
www.damtp.cam.ac.uk/user/majid
\end{center}

\begin{center}
June 1997
\end{center}
\vspace{10pt}
\begin{quote}\baselineskip 13pt
\noindent{\bf ABSTRACT}
We show that noncommutative differential forms on $k[x]$, $k$ a
field, are of the form $\Omega^1=k_\lambda[x]$ where
$k_\lambda\supseteq k$ is a field extension. We compute the case
$\C\supset \R$ explicitly, where $\Omega^1$ is 2-dimensional. We
study the induced quantum de Rahm complex, its cohomology and the
associated moduli space of flat connections.
\end{quote}
\baselineskip 23pt

\section{INTRODUCTION}

Let $A$ be an algebra, which we consider as playing the role of
`co-ordinates' in algebraic geometry, except that we do not require
the algebra to be commutative. The appropriate notion of cotangent
space or differential 1-forms in this case is\cite{Wor:dif}

1. $\Omega^1$ an $A$-bimodule

2. $\extd:A\to \Omega^1$ a linear map obeying the Leibniz rule
$\extd(ab)=a\extd b+(\extd a)b$ for all $a,b\in A$.

3. The map $A\tens A\to \Omega^1$, $a\tens b\mapsto a\extd b$ is
surjective.

When $A$ has a Hopf algebra structure with coproduct $\Delta:A\to
A\tens A$ and counit $\eps:A\to k$ ($k$ the ground field), we say
that $\Omega^1$ is {\em bicovariant} if

4. $\Omega^1$ is a bicomodule with coactions $\Delta_L:\Omega^1\to
A\tens\Omega^1,\Delta_R:\Omega^1\to \Omega^1\tens A$ bimodule maps (with
the tensor product bimodule structure on the target spaces, where
$A$ is a bimodule by left and right multiplication).

5. $\extd$ is a bicomodule map with the left and right regular
coactions on $A$ provided by $\Delta$.

A morphism of calculi means a bimodule and bicomodule map forming a
commuting triangle with the respective $\extd$ maps. One
says\cite{Ma:cla} that a calculus is {\em coirreducible} if it has
no proper quotients.

The main difference is that, in usual algebraic geometry, the
multiplication of forms $\Omega^1$ by `functions' $A$ is the same
from the left or from the right. However, if $a\extd b=(\extd b)a$
then by axiom 2. we have $\extd(ab-ba)=0$, i.e. we cannot naturally
suppose this when $A$ is noncommutative. This possible
noncommutativity of forms and `functions' is the main
generalisation featuring in the above axioms. We say that a
differential calculus is noncommutative or `quantum'  if the left
and right multiplication of forms by functions do not coincide. It
turns out that many geometrical constructions assume neither
commutativity of $A$ nor commutativity of the differential
calculus. See \cite{BrzMa:gau} for a theory of bundles with quantum
group fiber and the example of the $q$-monopole over a $q$-deformed
sphere. Moreover, there are natural prolongations to higher order
differential forms, i.e. the entire exterior algebra $\Omega^\cdot$
once $\Omega^1$ is specified.  Note that our approach is somewhat
different from \cite{Con:non}, where an entire $\Omega^\cdot$ on an
algebra is effectively specified via a `spectral triple'.

Bicovariant quantum differential calculi on strict quantum groups and group
algebras over $\C$ have recently been classified
in \cite{Ma:cla}. This class of structures is, however, interesting even
when $A$ is commutative, being a generalisation of usual concepts of
differential forms. In this paper we study the simplest case, where
$A=k[x]$, the polynomials over a field $k$, with its additive
coproduct and counit
\[ \Delta x=x\tens 1+1\tens x,\quad \eps x=0.\]
A complete classification of the bicovariant $\Omega^1$ in this
case is provided in Section~2. They turn out
to be of the form $\Omega^1=k_\lambda[x]$ where $k_\lambda\supset
k$ is a field extension. In Section~3 we study a concrete example
on $\R[x]$ associated to the field extension $\C\supset\R$.
Section~4 studies the exterior algebra, de Rahm cohomology and
elements of gauge theory in this setting. Section~5 concludes with some
further quantum geometric considerations.

\subsection*{Acknowledgements}
I would like to thank D. Solomon for useful discussions at the
start of the project.

\section{$\Omega^1$ and field extensions}

When $\Omega^1$ is required to be bicovariant, there is a standard
argument\cite{Wor:dif} that it must be of the form
\[ \Omega^1=\Omega_0\tens A,\quad \extd a=(\pi\tens\id)
(\Delta a-1\tens a)\] where $\Omega_0=\ker\eps/\CM$ with canonical
projection $\pi:\ker\eps\to\Omega_0$ and $\CM$ is a left ideal
contained in $\ker\eps$ and stable under the Hopf algebra adjoint
coaction $\Ad$. The right (co)module structures are those of $A$
alone by (co)multiplication. The left (co)module structures are the
tensor product of those on $\Omega_0$ as inherited from
$\ker\eps\subset A$ (where $A$ acts by left multiplication and
coacts by $\Ad$) and those on $A$ by (co)multiplication. We
recall that modules and comodules of a Hopf algebra $A$ have a
tensor product induced by the coproduct and product of $A$
respectively. Then bicovariant $\Omega^1$ are in 1-1 correspondence
with the $\Ad$-stable left ideals $\CM\subset\ker\eps$. When $A$ is
cocommutative the adjoint coaction $\Ad$ is trivial.

\begin{propos} When $A=k[x]$, the coirreducible bicovariant
$\Omega^1$ are in 1-1 correspondence with irreducible monic
polynomials $m\in k[x]$, and take the form $\Omega^1=k_\lambda[x]$
where $k_\lambda=k[\lambda]/\<m\>$ is the corresponding field
extension. The bimodule structures and $\extd$ differential are
\[ f(x)\cdot P(\lambda,x)=f(x+\lambda)P(\lambda,x),\quad
P(\lambda,x)\cdot f(x)=P(\lambda,x)f(x),\quad \extd
f(x)={f(x+\lambda)-f(x)\over\lambda}\] for all $f\in k[x],P\in
k_\lambda[x]$.
\end{propos}
\proof According to the above, bicovariant differential calculi
on $k[x]$ are  in 1-1 correspondence with ideals $\CM\subset
\ker\eps$. Here $\ker\eps=\<x\>$, the ideal generated by $x$ in
$k[x]$. Since $k[x]$ is a P.I.D., the ideal $\CM$ above is
generated by a polynomial. Since $\CM\subset\ker\eps$, this
polynomial is divisible by $x$, i.e. $\CM=\<xm\>$. Coirreducible
calculi correspond to $m$ irreducible and monic.

We identify the corresponding $\Omega_0=\<x\>/\<xm\>\isom
k[\lambda]/\<m\>=k_\lambda$ by $xf(x)\mapsto f(\lambda)$. Under
this identification, $\Omega^1=\Omega_0\tens k[x]\isom
k_\lambda[x]$. The action from the right is by the inclusion
$k[x]\subset k_\lambda[x]$. The action from the left is by
\[ f(x)\cdot x^m\tens x^n=f(x\tens 1+1\tens x)x^m\tens x^n.\]
as the tensor product action. Hence
$f(x)\cdot\lambda^{m-1}x^n=f(\lambda+x)\lambda^{m-1}x^n$ under our
identification. The quotient by $\<xm(x)\>$ or $\<m(\lambda)\>$ is
understood in these expressions.

We compute $\extd f=f(x\tens 1+1\tens x)-1\tens f(x)$ modulo
$\<xm\>$ in the first tensor factor. Under our isomorphism this is
$f(\lambda+x)-f(x)/\lambda$ modulo $\<m(\lambda)\>$. Note that
$\extd x=x\tens 1$ modulo $\<xm\>$ become $\extd x=1\in
k_\lambda[x]$.

To see explicitly that the correspondence here is indeed 1-1,
suppose that $k_{\lambda_1}[x]\isom k_{\lambda_2}[x]$ as quantum
differential calculi associated to $m_1(\lambda_1)$ and
$m_2(\lambda_2)$. Since the isomorphism is in particular a right
module map under $k[x]$, it restricts to the identity on $k[x]$.
And since the isomorphism forms a commutative triangle with the
$\extd$ maps, it identifies ${f(x+\lambda_1)-f(x)\over\lambda_1}$
with ${f(x+\lambda_2)-f(x)\over\lambda_2}$ for all $f$. Taking
$f=x^2$ we have $2x+\lambda_1$ and $2x+\lambda_2$ identified, hence
$\lambda_1,\lambda_2$ identified. Similarly by induction
$\lambda_1^n$ and $\lambda_2^n$ are identified for all $n\ge0$. One
can also use the left module map property to conclude this. Hence
$m_1(\lambda_2)=0$ in $k_{\lambda_2}$. Hence $m_2$ divides $m_1$.
Since $m_1$ is monic and irreducible, we conclude that $m_1=m_2$ as
required. The converse direction is clear. \endproof

This generalises the observation in \cite{Ma:gosb} that
coirreducible bicovariant quantum differential calculi over $\C[x]$
are parametrized by $\lambda_0\in \C$ (say). Here
$m(\lambda)=\lambda-\lambda_0$ and $\pi(\lambda)=\lambda_0$. Hence,
in this case,
\eqn{Cdif}{\extd f=\extd x {f(x+\lambda_0)-f(x)\over \lambda_0}.}
The ratio on the right should be understood as the coefficient of
$\lambda_0$ in $f(x+\lambda_0)-f(x)$, i.e. we include the usual
differential calculus as the case $\lambda_0=0$. More generally, if
the extension is Galois, the roots $\lambda_i$ of $m$ are as many
as its degree and are primitive elements of $k_\lambda$, i.e.
$k_\lambda\isom k[\lambda_i]$ by setting $\lambda=\lambda_i$, for
each $i$. This gives us different ways of thinking of the
differentials in Proposition~2.1 concretely as finite differences,
all of them equivalent via the action of the Galois group of
$k_\lambda$ automorphisms that permute the $\lambda_i$.

It is easy to verify that $\Omega^1$ in Proposition~2.1 is
bicovariant under the left and right coactions
\eqn{DeltaLR}{ \Delta_RP(\lambda,x)=P(\lambda,x+y)\in k_\lambda[x]
\tens k[y],\quad \Delta_LP(\lambda,x)=P(\lambda,y+x)\in k[y]
\tens k_\lambda[x]}
induced by the coproduct $\Delta$, as it must be by construction.
Here the coacting copy of $A$ is denoted by $k[y]$. The space
$\Omega_0$ is the subspace of $\Omega^1$  invariant under the left
coaction $\Delta_L$, again by the general theory. Clearly, the
dimension of $\Omega_0$ over $k$, which is the dimension of the
quantum differential calculus, is the degree of $m$, the degree of
the associated field extension. The elements $\{1=\extd
x,\lambda,\cdots,\lambda^{\deg(m)-1}\}$ of $k_\lambda[x]$ are a
basis of right-invariant 1-forms.

\section{Quantum differentials for the complex extension of the reals}

In this section we consider in detail the case $k=\R$ and
$m(\lambda)=\lambda^2+1$. Then $k_\lambda=\C$. The space of
right-invariant 1-forms has basis
\[ \Omega_0=\{1=\extd x,\quad \lambda=\extd x^2-2(\extd x)x\}\]
We use the notations $\extd x$ and $\omega\equiv \extd x^2-2(\extd
x)x=x\extd x-(\extd x)x$ (by the Leibniz rule) for these two
1-forms in what follows.

\begin{lemma} The left part of the bimodule structure on $\Omega^1$
in this basis is given by
\[ x\cdot \extd x=(\extd x)x+\omega,\quad x\cdot\omega
=\omega x-\extd x\]
\end{lemma}
\proof The first equality is the definition of $\omega$ (given the
Leibniz rule). The second depends on the irreducible polynomial $m$
according to $x\cdot\omega=(x+\lambda)\lambda=\lambda x-1=\omega
x-\extd x$. \endproof

\begin{propos} The exterior differential is given by
\[ \extd f(x)=(\extd x) \Im f(x+i)+ \omega(f(x)-\Re f(x+i))\]
where $f\in \R[x]$ is   continued to $\C$ and $\Im,\Re$ denote
imaginary and real parts. The left and right multiplication of
forms by functions are related by
\[ f(x)\cdot\pmatrix{\extd x&\omega}=\pmatrix{\extd x&\omega}
\pmatrix{\Re &-\Im\cr \Im&\Re}f(x+\imath)\]
\end{propos}
\proof This follows directly as an example of Proposition~2.1
on writing $\lambda=\imath$. Here we provide a more conventional direct
proof based on the more conventional description in Lemma~3.1. First we
write Lemma~3.1 in matrix form
\[ x\cdot\pmatrix{\extd x&\omega}=\pmatrix{\extd x&\omega}
\pmatrix{x &-1\cr 1&x}\]
Then $f(x)\cdot\extd x=(\extd x)f(x+\Lambda)^1{}_1+\omega
f(x+\Lambda)^2{}_1$ where $\Lambda=\pmatrix{0&-1\cr 1&0}$ and the
numerical indices denote the matrix element. We regard $x$ for
these purposes as multiplied by the identity matrix. Similarly for
$f(x)\cdot\omega$.

Now, by induction on the Leibniz rule,
\align{\extd x^m\equad &&=x^{m-1}\extd x+ x^{m-2}(\extd x)x
+\cdots +(\extd x)x^{m-1}\\ &&=\pmatrix{\extd x&\omega}
(x^{m-1}+x^{m-2}(x+\Lambda)+\cdots+x(x+\Lambda)^{m-2}
+(x+\Lambda)^{m-1})\pmatrix{1\cr 0}\\ &&=\pmatrix{\extd
x&\omega}({(x+\Lambda)^m-x^m\over \Lambda})
\pmatrix{1\cr 0}=\pmatrix{\extd x&\omega}((x+\Lambda)^m-x^m)
\pmatrix{0\cr -1}.}
This provides the formula
\[\extd f=\pmatrix{\extd x&\omega}(f(x+\Lambda)-f(x))
\pmatrix{0\cr -1}\]
Since $\Lambda^2=-1$, we then identify the $\Lambda^0$ and
$\Lambda^1$ parts  of $f(x+\Lambda)-f(x)$ with the real and
imaginary parts of $f(x+i)-f(x)$ as stated. Similarly for
$f(x)\cdot\extd x$ and $f(x)\cdot\omega$.
\endproof

To gain further insight into this differential calculus it is
useful to embed it in the 1-parameter family corresponding to
$m(\lambda)=\lambda^2+q^2$, $q\in\R$. This is isomorphic to the
above $q=1$ case for all $q\ne 0$ and not irreducible for $q=0$. It
does, however, have an interesting limit as $q\to 0$. Briefly, the
relevant formulae are
\eqn{qdifa}{   x\cdot\pmatrix{\extd x&\omega}
=\pmatrix{\extd x&\omega}\pmatrix{x &-q^2\cr 1&x},}
resulting in
\eqn{qdifb}{ \extd f=q^{-2}\omega(f(x)-\Re f(x+\imath q))
+q^{-1}(\extd x)\Im f(x+\imath q).}
 This has a limit as $q\to 0$:
\eqn{jetdifa}{ x\omega=\omega x,\quad x\extd x=(\extd x)x+\omega,\quad
\extd f=\omega \h f''+(\extd x)f'}
in terms of the usual newtonian derivative  $f'$. This is the 2-jet
calculus in \cite{Ma:gosb} whereby up to second order derivatives
are viewed as `first order' with respect to the new calculus and an
appropriate `braided derivation'\cite{Ma:cla} rule. We see that this calculus,
although not coirreducible, arises naturally as a degenerate limit
of coirreducibles corresponding to the extension $\R\subset \C$.

\section{Quantum cohomology and gauge theory of field extensions}

In this section, we consider two natural prolongations of the
$\Omega^1(k[x])$ associated to a field extension to `exterior
algebras' $\Omega^n(k[x])$ of degree $n>1$. We compute the
first quantum cohomology for each prolongation in the case of the
extension $\R\subseteq\C$, and the associated gauge theory.

We recall first that a differential graded algebra $\Omega^\cdot$
over a unital algebra $A$ means a graded algebra with degree zero
part $A$ itself, and $\extd:\Omega^\cdot\to\Omega^\cdot$ which
increases the degree by 1 and obeys $\extd^2=0$ and the graded
Leibniz rule. In other words, $\Omega^\cdot$ has the algebraic
properties of an `exterior algebra' in DeRahm theory and one may
likewise compute its `quantum de Rahm cohomology'. Thus,
\eqn{cohom}{ H^1=\{\omega\in\Omega^1|\ \extd \omega
=0\}/\{\extd a|\ a\in A\}.}

Given $\Omega^1$, its maximal prolongation is defined as follows.
First of all, we recall that view of Axiom~3 above, we can write
$\Omega^1$ as a quotient of the universal calculus
$\Omega^1_U=\ker(\cdot:A\tens A\to A)$ by a sub-bimodule $\CN$.
Here $\Omega^1_U$ has the obvious bimodule structure from $A\tens
A$ and $\extd_U a=a\tens 1-1\tens a$. (Note that when $A$ is a Hopf
algebra then $A\tens A\isom A\tens A$ by $a\tens b\mapsto (\Delta
a)b$ restricts to $\Omega^1_U\isom
\ker\eps\tens A$ and $\CN\isom \CM\tens A$ giving the description used
in Section~2). Moreover, $\Omega^1_U$ is the degree 1 part of a
canonical $\Omega^\cdot_U$ (albeit with trivial quantum
cohomology). Here $\Omega^n_U\subset A^{\tens n+1}$ as elements in
the joint kernel of all product maps $\cdot_i$ multiplying the
$i,i+1$'th copies of $A$. This can also be identified with
$\Omega^n_U=\Omega^1_U\tens_A\cdots\tens_A\Omega^1_U$ in the
obvious way; see \cite{Con:non}\cite{BrzMa:gau}. Here
\eqn{univdif}{ \extd_U(a_0\tens a_1\cdots\tens a_n)
=\sum_{i=0}^{n+1}(-1)^{n+1-i}a_0\tens
\cdots \tens a_{i-1}\tens 1\tens a_{i}\tens \cdots\tens a_n}
One may check that $\extd_U\circ\extd_U=0$. The product $\wedge$ of
$\Omega_U^\cdot$ is given by multiplication between the two
adjacent copies of $A$. A general $\Omega^\cdot$ over $A$ is a
quotient of $\Omega^\cdot_U$ by a differential graded ideal (i.e.
an ideal stable under $\extd_U$). Without loss of generality we
assume that the degree 0 part of the ideal is zero. The degree 1
part is some subbimodule $\CN\subseteq\Omega^1_U$ and conversely,
given $\CN$ the maximal prolongation is provided by the
differential ideal generated by $\CN$. Its degree 2 part is
$\CF=\Omega^1_U\wedge
\CN+\CN\wedge\Omega^1_U+\extd_U\CN$ and $\Omega^2=\Omega^2_U/\CF$.

\begin{lemma} For the field  extension $\R\subset\C$, the maximal
$\Omega^2$ is generated as an $\R[x]$-module by the two forms
$\extd x\wedge\extd x$ and $\extd x\wedge\omega$. Moreover,
\[ \extd \omega=2\extd x\wedge\extd x=2\omega\wedge\omega,
\quad \extd x\wedge\omega=-\omega\wedge\extd x.\]
\end{lemma}
\proof The subbimodule $\CN$ in our case is generated by
$x\omega-\omega x+\extd x$ where $\omega$ is defined as above. Now,
$\extd\omega=\extd(x\extd x-(\extd x)x)=2\extd x\wedge\extd x$ from
the definition of $\omega$ and the graded Leibniz rule and
$\extd^2=0$. Hence the subbimodule $\CF$ is generated by
$\Omega^1_U\wedge\CN,\CN\wedge\Omega^1_U$ and $\extd
x\wedge\omega+\omega\wedge
\extd x+2x\extd x\wedge\extd x- 2(\extd x\wedge\extd x)x$. From
Lemma~3.1 we have $x\extd x\wedge\extd x=(\extd x)\wedge x\extd
x+\omega\wedge\extd x=(\extd x\wedge\extd x)x+\extd
x\wedge\omega+\omega\wedge\extd x$ up to terms in
$\Omega^1_U\wedge\CN,\CN\wedge\Omega^1_U$. Therefore, $\CF$ is
generated by these and $\omega\wedge\extd x+\extd x\wedge\omega$.

Finally, from the definition of $\omega$, the relations in
$\Omega^1$ and $\CN$, we have
\align{\omega\wedge\omega\equad&&=(x\extd x-(\extd
x)x)\wedge\omega=-x\omega\wedge\extd x-(\extd x)x\wedge\omega\\
&&=\extd x\wedge\extd x-\omega\wedge x\extd x-(\extd
x)x\wedge\omega\\ &&=\extd x\wedge\extd x-\omega\wedge
\omega-\omega\wedge\extd x x+\extd x\wedge\extd x-(\extd x)\wedge
\omega x=2\extd x\wedge\extd x-\omega\wedge\omega}
which gives the stated description of $\Omega^2$ as a quotient of
the tensor square over $\R[x]$ of $\Omega^1$.
\endproof

\begin{propos} With the maximal $\Omega^2$, the quantum de
Rahm cohomology $H^1$ associated to $\R\subset\C$ vanishes.
\end{propos}
\proof Suppose $\extd((\extd x)f+\omega g)=0$ i.e.
$-\extd x\wedge\extd f+2(\extd x\wedge\extd x)g-\omega\wedge\extd
g=0$. Put in the form of $\extd f$ and $\extd g$ from
Proposition~3.2 and we see this is equivalent to
\[ \Im g(x+i)=f(x)-\Re f(x+i),\quad \Re g(x+i)+g(x)=\Im f(x+i)\]
which can be combined into the single equation
\eqn{cocyCR}{ f(x+i)-f(x)=i(g(x)+g(x+i)).}
We now show that such $f,g$ are necessarily of the form
\[ f=\Im h(x+i),\quad g=h(x)-\Re h(x+i)\]
for some $h(x)$. Note first that if $(f,g)$ obey (\ref{cocyCR}) and
without loss of generality $f=nx^{n-1}+$ lower degree, say, then
\[   g={n(n-1)\over 2}x^{n-2}+{\rm lower\ degree}.\]
Indeed,  writing $g=\mu x^p+$ lower degree, the second half of
(\ref{cocyCR}) implies that $\mu x^p+\mu\Re
(x+\imath)^p+\cdots=2\mu
x^p+\cdots=n\Im(x+\imath)^{n-1}-nx^{n-1}+\cdots=n(n-1)x^{n-2}+\dots$.
Equating leading terms gives $p=n-2$ and $2\mu=n(n-1)$.

Now let
\[ f_n=\Im (x+\imath)^n,\quad g_n=x^n-\Re (x+\imath)^n=x^n
-(x+\imath)^n+\imath f_n\]
for $n>0$. Note that the leading term of $f_n$ is $nx^{n-1}$ and
the leading term of $g_n$ is ${n(n-1)\over 2}x^{n-2}$. Hence
\[ f=f_n+\bar f,\quad g=g_n+\bar g\]
defines two polynomials $\bar f,\bar g$ of lower degree. Now since
$(f_n,g_n)$ are the components of the differential of $x^n$, and
since $\extd^2=0$, we know that they obey (\ref{cocyCR}). Hence
$(\bar f,\bar g)$ obeys (\ref{cocyCR}) and has lower degree.

Therefore we have a proof by induction. The case where $n=2$ is
easily seen to be true. I.e. if $f=2x+\mu$ then (\ref{cocyCR})
implies as above that $g=1$. Then indeed $f=\Im(
(x+\imath)^2+\mu(x+\imath))$ and $x^2+\mu x-\Re ( (x+\imath)^2+\mu
(x+\imath))=1=g$ as required. In terms of differential forms, the
assertion is that if $(\extd x)(2x+\mu)+\omega h(x)$ is closed then
$h(x)=1$ and the form is $\extd(x^2+\mu x)$. This may also be
verified directly from the relations in Lemma~4.1.
\endproof

Next we consider a natural quotient of the above prolongation which
always exists when $A$ is a Hopf algebra and $\Omega^1$ is
bicovariant. In this case $\Omega^1=\Omega_0\tens A$ as explained
in Section~2, and $\Omega^\cdot$ is defined in such a way that the
invariant differential forms `braided-anticommute' where the
braiding is the one associated to the quantum double of
$A$\cite{Wor:dif}. Fortunately, in our case where $A$ is
commutative and cocommutative, the quantum double braiding is the
trivial flip map (the usual transposition). Hence in this case we
have simply $\Omega^n=\Lambda^n\Omega_0\tens A$, where $\Lambda^n$
denotes the usual exterior algebra of the vector space $\Omega_0$.
We call this the {\em skew exterior algebra}.

\begin{propos} The skew $\Omega^2$ in the case $\R\subset\C$  is
1-dimensional with basis $\extd x\wedge\omega$ (i.e. as in
Lemma~3.1 with the additional relations $(\extd x)^2=0=\omega^2$).
The first quantum cohomology in this case is $H^1=\R\omega$, i.e.
1-dimensional and spanned by $\omega$.
\end{propos}
\proof This time $\extd((\extd x)f+\omega g)=0$ and $\extd f,\extd g$
from Proposition~3.2 implies only that
\eqn{cocywor}{ \Im g(x+\imath)=f(x)-\Re f(x+\imath)}
as the coefficient of $\extd x\wedge\omega$. (The first half of
(\ref{cocyCR}) no applies since $\extd x\wedge\extd x=0$.) This
equation still implies that if $f=nx^{n-1}+$ lower degree and $n>2$
then $g={n(n-1)\over 2}x^{n-2}+$ lower degree as before. Indeed, if
$g=\mu x^p+\cdots$ then it says $\mu p
x^{p-1}+\cdots=n{(n-1)(n-2)\over 2}+\cdots$. This is weaker than
before because it does not fix $\mu$ when $n=2$. We proceed as
before by writing $f=f_n+\bar f$, $g=g_n+\bar g$ so that $\bar
f,\bar g$ obey (\ref{cocywor}) and have lower degree. In this way
we obtain (without loss of generality by scaling $f,g$ suitably)
$f=F+2x+\mu$ and $g=G+\tau$ where $(\extd x)F+\omega G=\extd h$ for
some $h$. Adding $f_2+\mu f_1$ and $g_2=1$ (here $g_1=0$) to $F,G$,
we have $(\extd x)f+\omega g=(1-\tau)\omega+\extd h'$ for
$h'=h+x^2+\mu x$. Hence $H^1=\R \omega$. Indeed $\extd\omega=0$ for
this choice of $\Omega^2$ but $\omega$ is not exact.
\endproof

Finally,  associated to any $\Omega^\cdot$ over a unital algebra
$A$ one has further `quantum geometrical' constructions, such as
gauge theory. In its simplest form we consider a gauge field as any
$\alpha\in\Omega^1$ and a gauge transform as an any invertible
$\gamma\in A$. The group of gauge transforms acts on the set of
$\alpha$ by
\eqn{gaugetran}{ \alpha^\gamma=\gamma^{-1}\alpha\gamma
+\gamma^{-1}\extd\gamma}
The fundamental lemma of gauge theory is that the curvature
\eqn{gaugeF}{ F(\alpha)=\extd\alpha+\alpha\wedge\alpha\in\Omega^2}
is covariant in the sense
$F(\alpha^\gamma)=\gamma^{-1}F(\alpha)\gamma$. Moreover, one can
consider sections $\psi\in A$ and a covariant derivative $\nabla
\psi=\extd\psi+\alpha\psi\in \Omega^1$. One has an action of the group of
gauge transformations by $\psi^\gamma=\gamma^{-1}\psi$ and
$\nabla^\gamma\psi^\gamma=(\nabla\psi)^\gamma$. These facts require
only that $\Omega^1,\Omega^2$ obey the natural axioms as part of a
differential graded algebra, see\cite{Ma:gosb}. Note that when
$\Omega^1$ is `quantum', the nonlinearity in $F$ does not
necessarily collapse even though the `structure group' here is
trivial, i.e. one has many of the features of nonAbelian gauge
theory. One may also consider $\alpha$ with values in some other
algebra.

In our present setting where $A=k[x]$, only $1$ will be invertible
as a polynomial. One may enlarge $A$ and our constructions above to
handle this. Alternatively, instead of the `finite' gauge
transformations $\gamma$ one can consider only `infinitesimal'
ones. Here an infinitesimal gauge transformation means $\theta\in
k[x]$ acting by
\eqn{gaugeinf}{ \alpha^\theta=\alpha+\extd\theta+\alpha
\theta-\theta\alpha,\quad
F(\alpha^\theta)=F(\alpha)+F(\alpha)\theta-\theta F(\alpha)}
 to lowest order in $\theta$. This can be stated more formally as
a vector field associated to each $\theta$ on the space of
connections, etc., in the usual way. The covariant derivative
$\nabla=\extd+\alpha\wedge$ is covariant to lowest order under
$\psi^\theta=\psi-\theta\psi$. By the same methods as in
\cite{Ma:gosb} one may check that any $\Omega^1,\Omega^2$ which are
part of an exterior algebra will do for these features of gauge
theory. Covariance of the curvature means that the vector fields
associated to $\theta$ restrict to vector fields on the space of
flat connections. They may not, however, restrict to only the
algebraic (i.e. polynomial) part.

\begin{propos} For the extension $\R\subset\C$ and the maximal
prolongation $\Omega^2$ we write
 $\alpha=(\extd x)a+\omega b$ and
$F(\alpha)=(\extd x)^2F_0+\extd x\wedge\omega F_1$, say, then
\[ F_0+\imath F_1=(a(x)+\imath (b(x)+1))(a(x+\imath)
-\imath (b(x+\imath)+1))-1\]
and the infinitesimal gauge transformations are
\[ (a(x)+\imath (b(x)+1))\mapsto (a(x)+\imath (b(x)+1))(1+\theta(x)
-\theta(x+\imath)).\]
The algebraic part of the space of flat connections is a circle
\[ {\rm Flat}=\{\extd x s+ \omega t|\ s,t\in\R,\quad s^2+ (t+1)^2=1\}
=S^1\subset \C.\]
Here $s+\imath t\in \C\subset \C[x]\subset \Omega^1$ is a circle of unit
radius centered at $-\imath$. The action of $\theta(x)=x$ is a unit vector
field along the circle, so that the algebraic moduli space is the class of
the zero connection.
\end{propos}
\proof We use  Proposition~3.2 to compute $\extd\alpha+\alpha\wedge\alpha$
in $\Omega^2$. We then use Lemma~4.1 and collect the
coefficients of $\extd x\wedge\extd x$ and $\extd x\wedge\omega$ as
\cmath{F_0=(-\Im a(x+\imath)+\Re b(x+\imath))(1+b(x))+b(x)
+(\Re a(x+\imath)+\Im b(x+\imath))a(x)\\
\quad F_1=(\Re a(x+\imath)+\Im b(x+\imath))(1+b(x))-a(x)
+(\Im a(x+\imath)-\Re b(x+\imath))a(x).}
Likewise from Proposition~3.2, the action of infinitesimal gauge
transformation
$\theta\in \R[x]$ is
\cmath{ a\mapsto a(x)(1+\theta(x))+\Im \theta(x+\imath)(1+b(x))
-\Re\theta(x+\imath) a(x)\\
b\mapsto b(x)(1+\theta(x))+\theta -\Re\theta(x+\imath)(1+b(x))-\Im
\theta(x+\imath)a(x).}
We can then combine these expressions into the expressions shown
for $F_0+\imath F_1$ and $a+\imath b$. Note that $\alpha=\extd x
a+\omega b
=a+\imath b$ in the identification of Proposition~2.1, and similarly
$F(\alpha)=\extd x\wedge (F_0+\imath F_1)$ by an extension of this
identification.

Next we compute the algebraic part of the space of flat
connections. Suppose that
\[ \alpha=a+\imath b=s x^n+\cdots + \imath (t x^m+\cdots),\quad s,t\ne0\]
are the leading terms for the real and imaginary parts. Here $n,m\ge 0$.
Then
\align{F(\alpha)\equad&&=-s nx^{n-1}-\imath s{n(n-1)\over 2}
x^{n-2}+2 t x^m+\imath tmx^{m-1}\\
&&+(sx^n+\imath tx^m)(s x^n - \imath t x^m+\imath s
nx^{n-1}+tmx^{m-1}-s{n(n-1)\over 2}x^{n-2}+\imath t{m(m-1)\over
2}x^{m-2})+\cdots\\ &&=- s nx^{n-1}-\imath s{n(n-1)\over
2}x^{n-2}+2 t x^m+\imath tmx^{m-1}+s^2x^{2n}+t^2 x^{2m}\\ &&+\imath
s^2n x^{2n-1}+\imath t^2 m x^{2m-1}+ st(m-n)x^{m+n-1}+{\imath\over
2} st(m(m-1)-n(n-1))x^{m+n-2}+\cdots.} Now, since $n\ge 0$ the
$x^{n-1}$ and $\imath x^{n-2}$ terms can be dropped against the
$x^{2n}$ and $\imath x^{2n-1}$ terms respectively.

Suppose that $m\ge 1$. Then the $x^m$ and $\imath x^{m-1}$ terms
can likewise be dropped. If $m=n$ then
$(s^2+t^2)x^{2n}+\imath(s^2+t^2)nx^{2n-1}$ is dominant, in which
case $F=0$ would imply $s=0,t=0$. So this case is excluded under
our initial assumption. If $m>n$ then $t^2x^{2m}+\imath t^2
mx^{2m-1}$ is dominant, in which case $t=0$. Likewise $m<n$ would
imply $s=0$.

Hence $m=0$ for a flat connection under our assumption $s,t\ne 0$.
In this case, if $n\ge 1$ then  $ s^2 x^{2n} +\imath s^2nx^{2n-1}$
is dominant and $F=0$ would imply $s=0$. Hence $n=0$ as well for a
flat connection.

It remains to consider the simpler cases where $t=0$ or $s=0$ in
our leading terms (i.e. real or imaginary $\alpha$). If $t=0$ and
$s\ne 0$ we similarly conclude that $n=0$ for a non-zero flat
connection. And if $s=0$ and $t\ne 0$  then $m=0$ for a non-zero
flat connection in the same way. Hence for an algebraic connection
of zero curvature, we are left with $\alpha=s+\imath t$ for
$s,t\in\R$. Then
\[ F(\extd x s+\omega t)=\extd(\extd x s+\omega t)
+(\extd x s+\omega t)\wedge
(\extd x s+\omega t)=(t^2+2t+s^2) \extd x\wedge\extd x\]
via Lemma~4.1, which tells us that $s^2+(t+1)^2=1$ for zero curvature.

For $\theta(x)=x\eps$, where $\eps\in\R$, we have the infinitesimal
gauge transform
$s+\imath(t+1)\mapsto (s+\imath(t+1))(1-\imath\eps)$
to lowest order in $\eps$. This is an infinitesimal rotation of $s+\imath t$
about $-\imath$. \endproof

Although we can consider only infinitesimal gauge transformations
in our present algebraic setup, it is clear that the exponentiation
of the infinitesimal gauge transformations associated to
$\theta(x)=x\eps$ rotate us around the stated $S^1$. Since this
$S^1$ passes through the origin, we see that all the algebraic zero
curvature solutions stated are connected in this way to the zero
connection by finite gauge transformations. Note also that
infinitesimal gauge transformations by  $\theta(x)=x^n\eps$, $n>1$
take us out of the space of algebraic zero curvature connections.
This tells us that additional zero curvature connections beyond
those in the proposition certainly exist in a suitable context,
just not as polynomials. For example, the formal exponentiation of
the gauge transform by $\theta(x)=x^2\eps$ of the $\alpha=-2\imath$
solution is
\[ \alpha  =-\imath (1+e^{\tau(1-2x\imath)}),\quad \tau\in\R.\]
It corresponds to the gauge transformation of $\alpha=-2\imath$ by
$\gamma(x)=e^{\tau x^2}$, where (\ref{gaugetran}) for the $\R\subset\C$
calculus comes out as
\eqn{gaugeCR}{ \alpha^\gamma+\imath=(\alpha+\imath)
{\gamma(x)\over\gamma(x+\imath)}.}
Although we are not able to consider such finite gauge transformations and
exponentials in our polynomial setting, we see that the
infinitesimal gauge transforms do give us some information about
the entire space of solutions.

\begin{propos} For the extension $\R\subset\C$ and the skew
prolongation $\Omega^2$ we write
 $\alpha=(\extd x)a+\omega b$ as above and
$F(\alpha)=\extd x\wedge\omega F_1$, say. Then
\cmath{F_1=(\Re a(x+\imath)+\Im b(x+\imath))(1+b(x))-a(x)
+(\Im a(x+\imath)-\Re b(x+\imath))a(x)} and the infinitesimal gauge
transformations by $\theta$ as in the preceding proposition.  The
algebraic part of the space of flat connections is the complex
plane
\[ {\rm Flat}=\{\extd s+\omega t|\ s,t\in\R\}=\C\]
where $s+\imath t\in\C\subset\C[x]=\Omega^1$. The algebraic moduli
space of flat connections modulo gauge transformations is the
half-line $\R_+$.
\end{propos}
\proof  We take the same form for $\alpha$ with leading coefficients
$s,t$ as in the
preceding proof. This time, however, the zero curvature condition
is only
half of the preceding one. Indeed,
\align{F_1\equad&&=-s{n(n-1)\over 2}x^{n-2}+tmx^{m-1}
+tx^m(-s{n(n-1)\over 2}x^{n-2} +tmx^{m-1})\\
&&+sx^n(snx^{n-1}+t{m(m-1)\over 2}x^{m-2})+\cdots} for the leading
terms after cancellations. We used the same expression for $F_1$ as
the coefficient of $\extd x\wedge\omega$ in the
preceding proof. We drop $x^{n-2}$ against $x^{2n-1}$ and, assuming
$m\ge 1$ we drop $x^{m-1}$ as well. If $m=n$ we drop $x^{m+n-2}$
and the dominant term is $(s^2+t^2)nx^{2n-1}$, which would imply
$s=t=0$ for a flat connection. If $m>n$ the dominant term is
$t^2x^{2m-1}$ which would imply $t=0$. If $m<n$ the dominant term
is $s^2x^{2n-1}$ which would imply $s=0$. Hence $m=0$. Hence the
dominant term is $s^2nx^{2n-1}$ which would imply $s=0$ if $n\ge
1$. Hence $n=0$ as well. Finally, if we consider the similar form
of $\alpha$ with $s=0$, the leading term for $m\ge 1$ would be
$t^2x^{2m-1}$ and imply $\alpha=0$, so $m=0$ in this case for a
nonzero flat connection. If we consider $\alpha$ with $t=0$ then
the leading term is $s^2nx^{2n-1}$ as before, which would imply
$n=0$. These are similar arguments to those in the preceding proof
but relying now only on the imaginary part of the curvature. We
deduce that an algebraic flat connection is of the form
$\alpha=\extd x s+\omega t$. This time, however, $F(\alpha)=0$ for
all $s,t\in\R$ since $\extd x\wedge\extd x=0$ in the skew
prolongation.

Infinitesimal gauge transformations are computed as before without
change. Hence the ones of the form $\theta(x)=x\eps$ rotate about
$-\imath$ in the $s+\imath t$ plane. The orbits are circles of
constant radius $s^2+(t+1)^2\in\R_+$. The different orbits are
however inequivalent at least by such $\theta$. On the other hand,
higher degree $\theta$ take us out of the class of polynomial
connections. Hence the algebraic part of the moduli space of flat
connections is $\R_+$.
\endproof

Finally, the cohomology and moduli spaces in the maximal and skew
prolongations are much more easily computed in the simpler 2-jet
calculus resulting from the degenerate $q\to 0$ limit of the
parametrized version of the $\R\subset \C$ extension. We first
compute the maximal prolongation as having relations
\eqn{qdifc}{\omega\wedge\omega=q^2\extd x\wedge\extd x,\quad
\extd \omega=2\extd x\wedge\extd x,
\quad \extd x\wedge\omega=-\omega\wedge\extd x.}
The proof is entirely similar to that of Lemma~4.1 (and equivalent
to it after a rescaling), so we omit it. The degenerate limit is
therefore
\eqn{degdifb}{\omega\wedge\omega=0,\quad \extd\omega
=2\extd x\wedge\extd x,
\quad \extd x\wedge\omega=-\omega\wedge\extd x.}
The skew prolongation has the additional relation $\extd
x\wedge\extd x=0$.

\begin{propos} The quantum cohomology for the 2-jet calculus is $H^1=0$
in the maximal prolongation and $H^1=\R\omega$ in the skew
prolongation.
\end{propos}
\proof Here $\extd((\extd x)f+\omega g)=0$ implies
\[ \h f'=g,\quad \h f''=g'\]
Letting $h$ be such that $h'=f$, we have $(\extd x)f+\omega g=\extd
h$, so that $H^1$ is trivial. For the skew prolongation we have
only $\h f''=g'$, which implies $\h f'=g-\mu$ where $\mu\in\R$.
Choosing $h$ such that $f=h'$, we have $(\extd x)f+\omega g=\extd
h+\mu\omega$, so that $H^1=\R\omega$.
\endproof

Gauge theory   in the $q\to 0$ limit  is described in
\cite{Ma:gosb} and we now compute the moduli space of flat
connections in this case.

\begin{propos} Writing $\alpha=(\extd x)a+\omega b$, the curvature in the
2-jet calculus with the maximal prolongation is
\[ F(\alpha)=\extd x\wedge\extd x(2b-a'+a^2)
+\extd x\wedge\omega (b'-\h a''+a'a)\]
and is invariant under the gauge transformation
\[ a\mapsto
a+\theta',\quad b\mapsto b-a\theta'+\h(\theta''-(\theta')^2)\] by
$\theta\in\R[x]$. The moduli space of flat connections in the
maximal prolongation is trivial and in the skew prolongation is
$\R$, with flat connections gauge equivalent to $\alpha=\mu\omega$
for unique $\mu\in \R$.
\end{propos}
\proof We compute $F(\alpha)$ using the relations in $\Omega^2$ and
the commutation rules for $\Omega^1$ at the end of Section~2. The
gauge transformation is likewise the infinitesimal gauge
transformations as above but computed for this calculus, and
corrected by the $-\h(\theta')^2$ to make, in the present case, an
exact gauge symmetry of the curvature (not only to lowest order in
$\theta$)). These formulae are obtained by formally writing
$\gamma=e^\theta$ in the finite gauge transformation formulae
computed for the 2-jet calculus in \cite{Ma:gosb}; in our present
case the result involves only polynomials in derivatives of
$\theta$, i.e. makes sense in terms of $\theta$ at our algebraic
level. One then verifies directly at this level that
$F(\alpha^\theta)=F(\alpha)$.

The zero curvature condition in the maximal prolongation is
therefore
\[b=\h(a'-a^2).\]
If this is the case then choose $\theta$ such that $\theta'=-a$.
This gauge transforms $a\mapsto 0$. On the other hand, $b\mapsto
b-a(-a)+\h(-a'-a^2)=b+\h (a^2 -a')=0$ as well. Hence every flat
connection is gauge equivalent to the zero one. By contrast, in the
skew prolongation, the zero curvature equation is
\[ b'=\h a''-a'a\]
which means $b=\h(a'-a^2)+\mu$ for some constant $\mu\in\R$. Making
the same gauge transformation as before now sends $a\mapsto 0$ and
$b\mapsto\mu$. Any further gauge transformation preserving $a=0$
would require $\theta'=0$, which would therefore not change the $b$
component, i.e. the different $\mu$ cannot be related by any
further gauge transformation. Hence the moduli space is $\R$ in the
skew prolongation.  \endproof

\section{Concluding Remarks}

We conclude the paper with two miscellaneous pieces of general
theory, demonstrated for our particular quantum exterior algebras.
First, by \cite{Brz:rem}, the Woronowicz $\Omega^\cdot$ (which in
our case means the skew prolongation) is always a $\Z_2$-graded
Hopf algebra with coproduct extended by $\Delta=\Delta_L+\Delta_R$
on $\Omega^1$. The same applies in general to the maximal
prolongation, which again gives a $\Z_2$-graded Hopf algebra. From
(\ref{DeltaLR}), we know (for any field extension) that
$\Delta_L\lambda^n=1\tens\lambda^n$ and
$\Delta_R\lambda^n=\lambda^n\tens 1$ (i.e. $\Omega_0=k_\lambda$ is
left and right invariant). Hence the coproduct structure is with
the basis of $\Omega_0$ primitive, and the original coproduct of
$k[x]$.

For example, for the extension $\R\subset\C$ we have  the maximal
prolongation $\Omega^\cdot$ as the $\Z_2$-graded Hopf algebra
generated over $\R$ by $x$ of degree zero and $\theta\equiv\extd
x,\omega$ of degree 1, and the relations and coproduct
\ceqn{superHopf}{ x\theta-\theta x=\omega,\quad x\omega-\omega x
=-\theta,\quad \omega\theta=-\theta\omega,
\quad \theta^2=\omega^2\\
 \Delta x=x\tens 1+1\tens x,\quad
\Delta\theta=\theta\tens 1+1\tens\theta,
\quad \Delta\omega=\omega\tens1+1\tens\omega.}
The skew prolongation is the quotient of this by the additional
relations $\theta^2=0$.

Finally, we consider what should be the notion of `differentiable'
map $k[x]\to k[x]$ where the source and target are considered with
differential calculi defined by $m_1,m_2$ respectively. A full
analysis of the dependence of the above quantum geometric
constructions on the choice of $m$ will be developed elsewhere, but
one may conjecture that at least some `geometric' invariants
obtained from constructions of this type will be invariants of the
field extension; i.e. if $m_1,m_2$ give isomorphic field extensions
then some of the invariants should coincide. This is a long term
goal suggested by the above results, and would have applications in
number theory (where the question of which monic polynomials gives
equivalent extensions is poorly understood for many fields $k$).
The analysis of which maps $k[x]\to k[x]$ are indeed differentiable
should is a first step in this geometric programme.

We recall that any $\Omega^1(A)$ over a unital algebra $A$ is a
quotient $\Omega^1_UA/\CN_A$ of the universal 1-forms
$\Omega^1_UA\subset A\tens A$. Any algebra map $\phi: A\to B$
(between unital algebras $A,B$) clearly induces a map
$\phi\tens\phi:\Omega^1_UA\to\Omega^1_UB$. Given this situation, we
say that $\phi$ is {differentiable} if $\phi\tens
\phi$ descends to a map  $\Omega^1(A)\to \Omega^1(B)$. If so, we
denote the map
by $\phi_*$ and note that it obeys the commutative diagram
\eqn{phistar}{ \matrix{A&{\buildrel \phi\over \longrightarrow}&B\cr
\extd \downarrow& &\downarrow\extd\cr
\Omega^1(A)&{\buildrel \phi_*\over\longrightarrow}&\Omega^1(B)}}
\medskip
\noindent since the universal $\extd_U$ for $A,B$ clearly obey this. The
condition for differentiability is that
$(\phi\tens\phi)(\CN_A)\subseteq\CN_B$.

\begin{propos} In the setting of Proposition~2.1, an algebra map
$\phi:k[x]\to k[x]$  defined by $\phi(x)=\Phi\in k[x]$ is
differentiable with respect to calculi defined by
$m_1(\lambda_1),m_2(\lambda_2)$ on the source and target
respectively iff
\[ \extd\Phi =0,\quad {\rm or}\quad  m_1(\Phi(\lambda_2+x)-\Phi(x))=0\]
in $k_{\lambda_2}[x]$. Then $\phi_*(P(\lambda_1,x))
=(\extd \Phi)P(\Phi(\lambda_2+x)
-\Phi(x),\Phi(x))$, where the product is in $k_{\lambda_2}[x]$.
\end{propos}
\proof We use the explicit isomorphism $\theta:\Omega^1_UA\isom
\ker\eps\tens A$ provided by
$\theta(a\tens b)=a\o\tens a\t b$ and $\theta^{-1}(a\tens
b)=a\o\tens (Sa\t)b$ where $S$ is the antipode and $\Delta
a=a\o\tens a\t$ (summation understood). In view of this, the map
$\phi\tens\phi$ becomes the map $\phi_*^U:\ker\eps\tens A\to
\ker\eps \tens A$ as given by
\[ \phi_*^U(a\tens b)=\theta(\phi(a\o)\tens \phi(Sa\t)\phi(b))
=\phi(a\o)\o\tens \phi(a\o)\t\phi(Sa\t)\phi(b)\]
for all $a\in\ker\eps$ and $b\in A$. In the present setting, this
becomes
\[ \phi_*^U(yg(y)\tens f(x))=(\Phi(y+x)-\Phi(x))g(\Phi(y+x)
-\Phi(x))f(\Phi(x))\]
for polynomials $f,g$ (we write $A\tens A=k[y,x]$). As in the proof
of Proposition~2.1, we further identify the source
$\ker\eps=k[\lambda_1]$ by $yg(y)\mapsto g(\lambda)$. We likewise
identify the target $\ker\eps=k[\lambda_2]$ in the similar say.
With these identifications understood, we have
\[ \phi_*^U(g(\lambda_1)\tens f(x))={\Phi(\lambda_2+x)
-\Phi(x)\over\lambda_2}
g(\Phi(\lambda_2+x)-\Phi(x))f(\Phi(x)).\] This map descends to the
quotients $k_{\lambda_1}=k[\lambda_1]/\<m_1\>$ and
$k_{\lambda_2}=k[\lambda_2]/\<m_2\>$ iff
\[ \phi_*^U(m_1(\lambda_1)\tens 1)=0\]
in $k_{\lambda_2}[x]$, i.e. iff
\[ (\extd \Phi)\, m_1(\Phi(\lambda_2+x)-\Phi(x))=0\]
in $k_{\lambda_2}[x]$, where we used the description of $\extd$ in
the target calculus from Proposition~2.1. This is the condition
stated. Of the two possibilities, the second is more interesting in
view of the form of $\phi_*$.
\endproof

For example, for the differential calculus associated to
$\R\subset\C$ in the source and target, the differentiability
condition is
\eqn{difCR}{ \Phi(x+\imath)-\Phi(x)=\cases{ 0\cr\pm \imath}}
which at the algebraic level means $\Phi(x)=\pm x+\mu$ or
$\Phi(x)=\mu$ for $\mu\in\R$. If we allow non-polynomials then
other possibilities, such as $\Phi(x)=e^{2\pi x}$, certainly open
up. By contrast, for the degenerate $2$-jet calculus in the source
and target, the differentiability condition is automatically
satisfied for all $\Phi\in\R[x]$. Here $m(\lambda)=\lambda^2$ is
not irreducible but one can use the same formulae (the calculus is
merely not coirreducible). Then
$\Phi(\lambda+x)-\Phi(x)=\lambda\Phi'$ and
$m(\lambda\Phi')=\lambda^2(\Phi')^2=0$ for all $\Phi$.


\end{document}